\documentclass[a4paper,12pt]{article}
\usepackage[T2A]{fontenc}
\usepackage[utf8]{inputenc} 
\usepackage{amssymb,amsfonts}
\usepackage{amsfonts,amsmath,amssymb,amscd,latexsym}
\usepackage{srcltx}

\linethickness{0.4pt}

\newcommand{\bi}{\bibitem}
\newcommand{\nb}{\newblock}

\newcommand{\be}[1]{\begin{equation}\label{#1}}
\newcommand{\ee}{\end{equation}}

\newcommand{\la}{\langle\,}
\newcommand{\ra}{\,\rangle}

\newcommand{\prf}{{\bf Proof.}\ }

\newcommand{\hgt}{\mathop{\rm ht}}

\newtheorem{thm}{\quad Theorem}
\newtheorem{lm}{\quad Lemma}
\newtheorem{cy}{\quad Corollary}

\newtheorem{prop}{\quad Proposition}
\newtheorem{que}{\quad Question}

\title{Systems of equations over the group ring of Thompson's group $F$}
\author{\vspace{2ex}
V. S. Guba\thanks{This work is supported by the Russian Foundation
for Basic Research, project no. 20-01-00465.}\\
Vologda State University,\\
15 Lenin Street,\\
Vologda\\
Russia\\
160600\\
E-mail: gubavs{@}vogu35.ru}
\date{}

\begin{document}

\maketitle

\begin{abstract}

Let $R=K[G]$ be a group ring of a group $G$ over a field $K$. It is known that if $G$ is amenable then $R$ satisfies the Ore condition: for any $a,b\in R$ there exist $u,v\in R$ such that $au=bv$, where $u\ne0$ or $v\ne0$. It is also true for amenable groups that a non-zero solution exists for any finite system of linear equations over $R$, where the number of unknowns exceeds the number of equations. Recently Bartholdi proved the converse. As a consequence of this theorem, Kielak proved that  R.\,Thompson's group $F$ is amenable if and only if it satisfies the Ore condition. The amenability problem for $F$ is a long-standing open question.

In this paper we prove that some equations or their systems have non-zero solutions in the group rings of $F$. We improve some results by Donnelly showing that there exist finite sets $Y\subset F$ with the property $|AY| < \frac43|Y|$, where $A=\{x_0,x_1,x_2\}$. This implies some result on the systems of equations. We show that for any element $b$ in the group ring of $F$, the equation $(1-x_0)u=bv$ has a non-zero solution. The corresponding fact for $1-x_1$ instead of $1-x_0$ remains open. We deduce that for any $m\ge1$ the system $(1-x_0)u_0=(1-x_1)u_1=\cdots=(1-x_m)u_m$ has nonzero solutions in the group ring of $F$. We also analyze the equation $(1-x_0)u=(1-x_1)v$ giving a precise explicit description of all its solutions in $K[F]$. This is important since to any group relation between $x_0$, $x_1$ in $F$ one can naturally assign such a solution. So this can help to estimate the number of relations of a given length between generators.

\end{abstract}

\section{Introduction}

Here we repeat some basic definitions and notation used in the paper. A more detailed information on the subject can be found in~\cite{Gu04,Gu21a,Gu21b}.

\subsection{Thompson's monoid and group}
\label{thmg}

Let $M$ be a monoid given by the following infinite presentation

\be{xinf}
\la x_0,x_1,x_2,\ldots\mid x_j{x_i}=x_ix_{j+1}\ (0\le i < j)\,\ra.
\ee

It is well-known that any word in these generators can be reduced to a word of the form $x_{i_1}\ldots x_{i_k}$, where $k\ge0$ and $0\le i_1\le\cdots\le i_k$. This is called a {\em normal form} of a monoid element.
\vspace{1ex}

For any $a,b\in M$ there exists a least common right multiple of $a$, $b$. Also it is easy to check that $M$ is a cancellative monoid. A classical Ore theorem states that for a cancellative monoid $M$ with common right multiples, there exists a natural embedding of $M$ into its group of quotients (see~\cite{Ljap} for details). Any element of this group belongs to $MM^{-1}$, and the group is given be the same presentation. We denote it by $F$.
\vspace{1ex}

This group was discovered by Richard J. Thompson in the 60s. Detais can be found in the survey \cite{CFP}; see also~\cite{BS,Bro,BG}. It is easy to see that $x_n=x_0^{-(n-1)}x_1x_0^{n-1}$ for any $n\ge2$, so the group is generated by $x_0$, $x_1$. It can be given by the following presentation with two defining relations

\be{x0-1}
\la x_0,x_1\mid x_1^{x_0^2}=x_1^{x_0x_1},x_1^{x_0^3}=x_1^{x_0^2x_1}\ra,
\ee
where $a^b\leftrightharpoons b^{-1}ab$. 

Elements of the monoid $M$ are called {\em positive elements}. Each element of $F$ can be uniquely represented by a group {\em normal form\/}, that is, an expression of the form $pq^{-1}$, where $p$, $q$ are normal forms of positive elements, 
and the following is true: if the word $pq^{-1}=\ldots x_i\ldots x_i^{-1}\ldots$ contains both $x_i$ and $x_i^{-1}$ for some $i\ge0$, then it also contains $x_{i+1}$ or $x_{i+1}^{-1}$.
\vspace{1ex}

Many authors use an equivalent definition of $F$ in terms of piecewise-linear functions. In fact, there are various presentations of $F$ in these terms. In our paper we are not going to use these functions at all so we do not list all conditions for them. This information is standard and can be found in may sources including~\cite{CFP,BS}. 
\vspace{1ex}

We will also need another finite presentation of $F$ in terms of symmetric generators $x_1$, $\bar{x}_1$, where $\bar{x}_1=x_1x_0^{-1}$. The symmetric system of generators $\{x_1,\bar{x}_1\}$ is sometimes more convenient since there exist an involutive automorphism of $F$ that takes $x_1$ to $\bar{x}_1$. 

Let $a\leftrightarrow b$ mean that $a$ and $b$ commute. Taking inverses of the above generators, we get the following presentaion of $F$:

\be{ab}
\la\alpha,\beta\mid\alpha^{\beta}\leftrightarrow\beta^{\alpha},\alpha^{\beta}\leftrightarrow\beta^{\alpha^2}\ra.
\ee
where $\alpha=x_1^{-1}$, $\beta=\bar{x}_1^{-1}=x_0x_1^{-1}$.

More details on this presentation can be found in~\cite{Gu21a}. In particular, for any $m,n\ge1$ it holds $\alpha^{\beta^m}\leftrightarrow\beta^{\alpha^n}$ as a consequence of the defining relations. 
\vspace{1ex}

It is known that $F$ is a {\em diagram group} over the simplest semigroup presentation ${\cal P}=\langle x=x^2\rangle$. See~\cite{GbS} for the theory of these groups. In our paper, we will concern diagram presentations for the elements of $F$, as well as marked binary forests for the same purposes. These techniques will be explained below.

\subsection{Graphs and automata}
\label{grr}

We will often work with graphs, including Cayley graphs of groups and also semigroup diagrams. By a {\em graph} we mean a (non-oriented) graph in the sense of Serre~\cite{Se80}. This means that each geometric edge consists of two mutually inverse directed edges. We do not describe it more formally since this concept is standard.
\vspace{1ex}

Let $G$ be a group generated by a set $A$. Its {\em right Cayley graph} $\Gamma_r={\cal C}(G;A)$ is defined as follows. The set of vertices is $G$; for any $a\in A^{\pm1}$ we put a directed edge $e=(g,a)$ from $g$ to $ga$. The letter $a$ is called the {\em label} of this edge. Labels are naturally extended to the set of paths in Cayley graphs.

In some cases we will prefer to work with {\em left} Cayley graphs $\Gamma_l={\cal C}(G;A)$. They are defined in a similar way. Here the set of vertices is also $G$. A directed egde $e=(a,g)$ with label $a$ goes from $ag$ to $g$. That is, going along the edge labelled by $a$, means cancelling $a$ on the left. Labels of paths are defined similarly for this case. It is worth noting that a  word $w$ in group generators will be the label of a path from the vertex $g\in G$ represented by $w$ to the identity element in case of left Cayley graphs.
\vspace{1ex}

The cardinality of a finite set $Y$ will be denoted by $|Y|$. 

Let $G$ be a group generated by $A$. For our needs we can assume that $A$ is always finite, $|A|=m$.  Let $\Gamma={\cal C}(G,A)$ be the Cayley graph of $G$, right or left. To any finite nonempty subset $Y\subset G$ we assign a subgraph in $\Gamma$ adding all edges connecting vertices of $Y$. So given a set $Y$, we will usually mean the corresponding subgraph. This is a labelled graph that we often call an {\em automaton}. For each $g\in Y$ we have exactly $2m$ directed edges in $\Gamma$ starting at $g$, where $a\in A^{\pm1}$. If the endpoint of such an edge with label $a$ belongs to $Y$, then we say that the vertex $g$ of our automaton $Y$ {\em accepts} $a$. For the case of right Cayley graphs this means $ga\in Y$, for the case of left Cayley graphs this means $a^{-1}g\in Y$. 

A vertex $g\in Y$ is called {\em internal} whenever it accepts all labels $a\in A^{\pm1}$. That is, the degree of $g$ in $Y$ equals $2m$. Otherwise we say that $g$ belongs to the {\em inner boundary} of $Y$ denoted by $\partial Y$. 

An edge $e$ is called {\em internal} whenever it connect two vertices of $Y$. If an edge $e$ connects a vertex of $Y$ with a vertex outside $Y$, then we call it {\em external}. The set of external edges form the {\em Cheeger boundary} of $Y$ denoted by $\partial_{\ast}Y$. 
\vspace{1ex}

By the {\em density} of a finite subgraph $Y$ we mean its average vertex degree. This concept was introduced in~\cite{Gu04}; see also~\cite{Gu21a}. It is denoted by $\delta(Y)$. A {\em Cheeger isoperimetric constant} of the subgraph $Y$ is the quotient $\iota_*(Y)=|\partial_{\ast}Y|/|Y|$. It follows directly from the definitions that $\delta(Y)+\iota_*(Y)=2m$. Indeed, each vertex $v$ has degree $2m$ in the Cayley graph $\Gamma$. This is the sum of the number of internal edges starting at $v$, which is $\deg_Y(v)$, and the number of external edges starting at $v$. Taking the sum over all $v\in Y$, we have $2m|Y|$, which is equal to $\sum\limits_v\deg_Y(v)+|\partial_{\ast}Y|$. Dividing by $|Y|$, we get the above equality.
\vspace{1ex}


\subsection{Semigroup diagrams}
\label{sgpd}

Let ${\cal P}$ be a semigroup presentation. Here it suffices to work with ${\cal P}=\la x\mid x=x^2\ra$. We will use 
non-spherical diagrams over ${\cal P}$ to represent elements of the group $F$. A more detailed description can be found in~\cite{Gu04}. Here we will mostly use diagrams that represent elements of the monoid $M$. Such objects are called {\em positive diagrams} over $x=x^2$. Let us give a brief illustration.

Given a normal form $pq^{-1}$ of an element in $F$, it is easy to draw the corresponding non-spherical diagram,
and vice versa. The following example illustrates the diagram that corresponds to the element $g=x_0^2x_1x_6x_3^{-1}x_0^{-2}$ represented by its normal form:

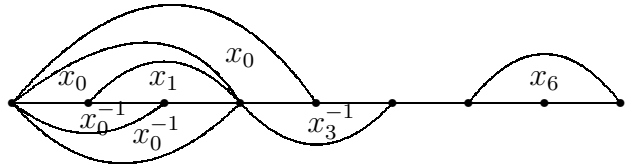
\begin{center}
	\begin{picture}(87.00,30.00)
	\put(6.00,9.00){\circle*{1.00}}
	\put(16.00,9.00){\circle*{1.00}}
	\put(26.00,9.00){\circle*{1.00}}
	\put(36.00,9.00){\circle*{1.00}}
	\put(46.00,9.00){\circle*{1.00}}
	\put(56.00,9.00){\circle*{1.00}}
	\put(66.00,9.00){\circle*{1.00}}
	\put(76.00,9.00){\circle*{1.00}}
	\put(86.00,9.00){\circle*{1.00}}
	\put(6.00,9.00){\line(1,0){80.00}}
	\bezier{132}(66.00,9.00)(76.00,22.00)(86.00,9.00)
	\bezier{120}(16.00,9.00)(26.00,20.00)(36.00,9.00)
	\bezier{176}(36.00,9.00)(26.00,25.00)(6.00,9.00)
	\bezier{264}(46.00,9.00)(28.00,35.00)(6.00,9.00)
	\bezier{120}(36.00,9.00)(46.00,-2.00)(56.00,9.00)
	\bezier{104}(6.00,9.00)(16.00,1.00)(26.00,9.00)
	\bezier{176}(6.00,9.00)(20.00,-7.00)(36.00,9.00)
	\put(36.00,15.00){\makebox(0,0)[cc]{$x_0$}}
	\put(14.00,12.00){\makebox(0,0)[cc]{$x_0$}}
	\put(26.00,12.00){\makebox(0,0)[cc]{$x_1$}}
	\put(76.00,12.00){\makebox(0,0)[cc]{$x_6$}}
	\put(25.00,5.00){\makebox(0,0)[cc]{$x_0^{-1}$}}
	\put(48.00,6.50){\makebox(0,0)[cc]{$x_3^{-1}$}}
	\put(18.00,7.00){\makebox(0,0)[cc]{$x_0^{-1}$}}
	\end{picture}
\end{center}

Here we assume that each edge of the diagram is labelled by a letter $x$. The horizontal path in the picture cuts the diagram into two parts, {\em positive} and {\em negative}. The positive part represents the element $x_0^2x_1x_6\in M$. The top path of this positive diagram has label $x^4$, the bottom path of it has label $x^8$. So we have a positive $(x^4,x^8)$-diagram over ${\cal P}$. The set of all positive diagrams with top label $x^m$ and bottom label $x^n$, where $m\le n$, will be denoted by ${\cal S}_{m,n}$. Given a positive $(x^m,x^n)$-diagram $\Delta_1$ and a positive $(x^n,x^k)$-diagram $\Delta_2$, one can {\em concatenate} them obtaining an $(x^m,x^k)$-diagram denoted by $\Delta_1\circ\Delta_2$. This operation is very natural: we identify the bottom path of $\Delta_1$ with the top path of $\Delta_2$. 

We will also say {\em product} instead of concatenation, and can multiply sets of diagrams in this way. Clearly, ${\cal S}_{m,n}{\cal S}_{n,k}={\cal S}_{m,k}$. All sets ${\cal S}_{m,n}$ are finite, and the cardinalities of them are given by numbers from the Catalan triangle.
\vspace{2ex}

Equivalently, positive elements of $F$ can be described in terms of rooted binary forests. One can use dual graphs for positive semigroup diagrams. However, we need presentation of elements in $F$ in terms of marked rooted binary forests from~\cite{BB05}. To be more formal, let us define a {\em rooted binary tree} by induction.
\vspace{0.5ex}

1) A dot $\cdot$ is a rooted binary tree.

2) If $T_1$, $T_2$ are rooted binary trees, then $(T_1\,\widehat\ \ T_2)$ is a rooted binary tree.

3) All rooted binary trees are constructed by the above rules.
\vspace{1ex}

Instead of formal expressions, we will use their geometric realizations. A dot will be regarded as a point. It coincides with the root of that tree. If $T=(T_1\,\widehat\ \ T_2)$, then we draw a {\em caret\/} for $\,\widehat\ \ $ as a union of two closed intervals $AB$ (goes left down) and $AC$ (goes right down). The point $A$ is the {\em root} of $T$. After that, we draw trees for $T_1$, $T_2$ and attach their roots to $B$, $C$ respectively in such a way that they have no intersection. It is standard that
for any $n\ge0$, the number of rooted binary trees with $n$ carets is equal to the $n$th Catalan number $c_n=\frac{(2n)!}{n!(n+1)!}$.
\vspace{1ex}

Each rooted binary tree has {\em leaves\/}. Formally, they are defined as follows: for the one-vertex tree
(which is called {\em trivial\/}), the only leaf coincides with the root. In case $T=(T_1\,\widehat\ \ T_2)$, the set of leaves
equals the union of the sets of leaves for $T_1$ and $T_2$. In this case the leaves are exactly vertices of degree
$1$.

We will also need the concept of a {\em height\/} of a rooted binary tree. For the trivial tree, its height equals
$0$. For $T=(T_1\,\widehat\ \ T_2)$, its height is $\hgt T=\max(\hgt T_1,\hgt T_2)+1$.

Now we define a {\em rooted binary forest\/} as a finite sequence of rooted binary trees $T_1$, ... , $T_m$,
where $m\ge1$. The leaves of it are the leaves of the trees. It is standard from combinatorics that the number
of rooted binary forests with $n$ leaves also equals $c_n$. The trees are enumerated from left to right and they
are drawn in the same way.

A {\em marked\/} (rooted binary) forest is a (rooted binary) forest where one of the trees is marked.

\subsection{Amenability}
\label{amen}
Recall that a group $G$ is called {\em amenable\/} whenever there exists a finitely additive normalized invariant mean on $G$. We are not going to list all well-known properties of amenable groups and non-amenable groups. We can just refer to one of surveys like~\cite{CGH} or Section 5.8 of~\cite{Sap14}.

The class of amenable groups is invariant under taking subgroups, quotient groups, group extensions, and directed unions of groups. It includes includes the class of all finite groups and all Abelian groups. The closure of that class under the above operations is the class EA of {\em elementary amenable\/} groups. A free group of rank $>1$ is not amenable. There are many useful
criteria for (non)amenability. Here we would like to mention Folner criterion from~\cite{Fol}. Let us formulate it as follows.
\vspace{1ex}

{\sl A group $G$ is amenable if and only if for any finite set $A\subseteq G$ and for any $\varepsilon > 0$, there exists a finite set $S$ such that $|AS| < (1+\varepsilon)|S|$.}
\vspace{1ex}

One can assume that $A$ contains the identity element. Also one can extend $A$ to the generating set of $G$ provided the group is finitely generated. In this case we see that the set $S$ from the above statement is almost invariant under (left) multiplication by elements in $A$.
\vspace{1ex}

It was proved in~\cite{BS} that the group $F$ has no free subgroups of rank $>1$. It is also known that $F$ is not elementary amenable (see~\cite{CGH}). However, the famous problem about amenability of $F$ remains open. The question whether $F$ is amenable was asked by Ross Geoghegan in 1979; see~\cite{Ger87}. There were many (unsucessful) attempts of various authors to solve this problem in both directions. However, to emphasize the difficulty of the question, we mention the paper~\cite{Moore13}, where it was shown that if $F$ is amenable, then Folner sets for it have a very fast growth. Besides, we would like to refer to the paper~\cite{BB05} where the authors obtained an estimate of the (Cheeger) isoperimetric constant of the group $F$ in its standard set of generators $\{x_0,x_1\}$. This estimate has not been improved so far.
\vspace{1ex}

\subsection{Ore condition and systems of equations}
\label{oca}

Tamari~\cite{Ta54} shows that if a group $G$ is amenable, then the group ring $R=K[G]$ satsfies Ore condition for
any field $K$. This means that for any $a,b\in R$ there exist $u,v\in R$ such that $au=bv$, where $u\ne0$ or $v\ne0$. That is, the equation $au=bv$ in the group ring has a non-zero solution.

This can be generalized as follows. Let we have a system of linear equations with coefficients in $R$, where the number of variables exceeds the number of equations:
$$
\left\{
\begin{array}{ccccccc}
a_{11}u_1&+  &\cdots  &+  &a_{1n}u_n  & = & 0 \\
\cdots&  &  \cdots & &\cdots  &  &  \\
a_{m1}u_1&+  &\cdots  &+  &a_{mn}u_n  &=  &0 
\end{array}
\right.
$$
Here $n > m$, $a_{ij}\in R$ for all $1\le i\le m$, $1\le j\le n$. We are interested in solutions $(u_1,...,u_n)\in R^n$. 
\vspace{1ex}

For amenable group $G$, this system always has a nonzero solution. This follows from cardinality arguments using the idea from~\cite{Ta54}. See~\cite{Gu21b} for details. 
\vspace{1ex}

In a recent paper~\cite{Ba19}, Bartholdi shows that the converse to the above statement is true. This gives a new criterion for amenabilty of groups. Although~\cite[Theorem 1.1]{Ba19} concerns the so-called GOE and MEP properties of automata (Gardens of Eden and Mutually Erasable Patterns), the proof of it allows one to extract the following statement.

\begin{prop}
\label{barth19}
{\rm(Bartholdi)}\ For any group $G$, the following two properties are equaivalent.

(i) $G$ is amenable

(ii) For any field $K$ and for any system of $m$ linear equations over $R=K[G]$ in $n > m$ variables, there exists a nonzero solution.
\end{prop}

In the Appendix to the same paper, Kielak shows that if the group ring $K[G]$ has no zero divisors, both properties are equivalent to the Ore condition. In particular, this holds for R.\,Thompson's group $F$. It is orderable, so there are no zero divisors in a group ring over a field. So we quote the following

\begin{prop}
\label{kielak}
{\rm(Kielak)}\ The group $F$ is amenable if and only if the group ring $K[F]$ over any field satisfies the Ore condition.
\end{prop}

The following elementary property of the group $F$ is well known.

\begin{lm}
\label{gig}
For any $g_1,...,g_n\in F$ there exist $g\in M$ such that $g_1g,\ldots,g_ng\in M$.
\end{lm}

It implies that equations and their systems over the group ring $K[F]$ can be reduced to the case of the monoid ring $K[M]$. It also  implies that any finite subset $Y\subset F$ can be translated into $M$ via a right multiplication by an element $g\in M$, that is, $Yg\subset M$. 

\section{Problem ${\cal Q}_3$}
\label{q3}

In~\cite{Gu21b} we introduced the following family of problems.
\vspace{1ex}

{\bf Problem} ${\cal Q}_k$: {\sl Given $k+1$ linear combinations of elements $x_0$, $x_1$, $x_2$, consider a system of $k$ equations with $k+1$ unknowns:
	$$
	(\alpha_0x_0+\beta_0x_1+\gamma_0x_2)u_0=(\alpha_1x_0+\beta_1x_1+\gamma_1x_2)u_1=\cdots=(\alpha_kx_0+\beta_kx_1+\gamma_kx_2)u_k.
	$$
	Find a nonzero solution of this system, where $u_0,u_1,\ldots,u_k\in K[M]$, or prove that it does not exist.}
\vspace{1ex}

We mentioned that Problem ${\cal Q}_2$ can be solved in positive since there exists a finite set $Y\subset F$ with the property $|AY| < \frac32|Y|$, where $A=\{x_0,x_1,x_2\}$. This can be done by cardinality reasons. Also we had announced a stronger fact. Namely, there exists a finite set $Y\subset F$ with the property $|AY| < \frac43|Y|$. This implies that Problem ${\cal Q}_3$ has a positive solution. In this Section we are going to give a proof. Notice that the size of the set $Y$ is very huge.
\vspace{1ex}

Let us recall the result of Donnelly from~\cite{Don14}. He shows that $F$ is non-amenable if and only if there exists $\varepsilon > 0$ such that for any finite set $Y\subset F$, one has $|AY|\ge(1+\varepsilon)|Y|$, where $A=\{x_0,x_1,x_2\}$ (see also~\cite{Don07}). For the set $Y$ here, one can assume without loss of generality that $Y$ is contained in ${\cal S}_{4,n}$ for some $n$. This gives some evidence that the amenability problem for $F$ has very close relationship with the family of Problems ${\cal Q}_k$. The case $k=4$ looks as a possible candidiate to a negative solution (that is, all solutions are zero). If true, this will imply that the constant $\varepsilon=\frac14$ fits into the above condition.
\vspace{1ex}

\begin{lm}
\label{xy}
Let $y_0=x_0x_2^{-1}$, $y_1=x_1x_2^{-1}$. Then the mapping $x_0\mapsto y_0$, $x_1\mapsto y_1$ induces a monomorphism of $F$ into itself.
\end{lm}

\prf Let $y_n=y_0^{-(1-n)}y_1y_0^{n-1}$ for $n\ge2$. Direct calculations show that $y_2=x_2x_0^{-1}x_1x_2^{-1}x_0x_2^{-1}=x_2^2x_3^{-1}x_2^{-1}$, $y_3=x_2x_0^{-1}x_2^2x_3^{-1}x_2^{-1}x_0x_2^{-1}=x_2x_3^2x_4^{-1}x_3^{-1}x_2^{-1}$, and $y_4=x_2x_0^{-1}x_2x_3^2x_4^{-1}x_3^{-1}x_2^{-1}x_0x_2^{-1}=x_2x_3x_4^2x_5^{-1}x_4^{-1}x_3^{-1}x_2^{-1}$. Now it follows that $y_2^{y_1}=x_2y_2^{x_1}x_2^{-1}=y_3$ and $y_3^{y_1}=x_2y_3^{x_1}x_2^{-1}=y_4$. So the defining relations of $F$ from~(\ref{x0-1}) hold via the replacements $x_i\mapsto y_i$ ($i\ge0$). Hence we have an endomorphism of $F$. Its image is non-Abelian since $y_2=y_0^{-1}y_1y_0\ne y_1$ as we can see from normal forms. According to~\cite[Theorem 4.13]{Bro}, the group $F$ has no proper non-Abelian homomorphic images. Therefore, we have a monomorphism. This completes the proof.
\vspace{1ex}

\begin{thm}
\label{43}
Let $A=\{x_0,x_1,x_2\}$. Then there exists a finite subset  $S\subset F$ satisfying $|AS| < \frac43|S|$.
\end{thm}

\prf We start from quite a general construction demonstrating self-similarity of $F$ in terms of its Cayley graphs.

Let $\Gamma={\cal C}(F,A)$ be the left Cayley graph of $F$ in generators $x_0$, $x_1$, $x_2$. For any path labelled by $x_0x_2^{-1}$ we add a new directed edge labelled by $y_0$ from the initial to the terminal vertex of this path. Similarly, for any path labelled by $x_1x_2^{-1}$ we add a new directed edge labelled by $y_1$ from the initial to the terminal vertex of this path. By Lemma~\ref{xy}, we get a left Cayley graph of $F$ in generators $y_0$, $y_1$ if we forget about all edges labelled by $x_i^{\pm1}$ ($0\le i\le2)$. Now let $\bar{y}_1=y_1y_0^{-1}=x_1x_0^{-1}=\bar{x}_1$. We add directed egdes labelled by $\bar{y}_1$ getting triangles like this:

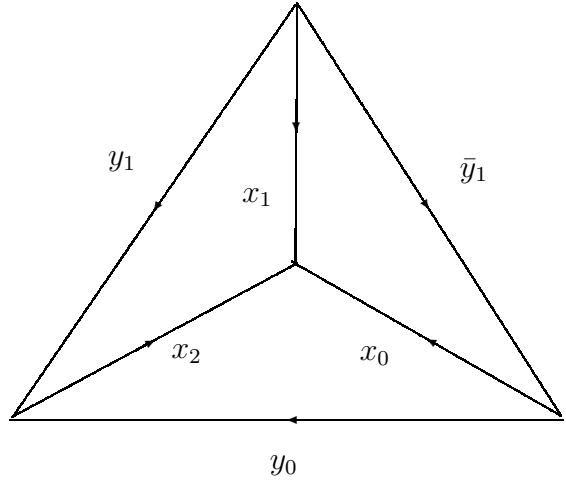
\begin{center}
\begin{picture}(76.5,69)(0,0)
\put(40.125,13.75){\vector(-1,0){.07}}\put(76.5,13.75){\line(-1,0){72.75}}
\put(22.625,41.375){\vector(-2,-3){.07}}\multiput(41.25,68.75)(-.03371040724,-.04954751131){1105}{\line(0,-1){.04954751131}}
\put(58.875,41.625){\vector(2,-3){.07}}\multiput(41.5,69)(.03373786408,-.05315533981){1030}{\line(0,-1){.05315533981}}
\put(58.5,24.5){\vector(-2,1){.07}}\multiput(76.25,14.25)(-.0583881579,.0337171053){608}{\line(-1,0){.0583881579}}
\put(22.875,24.375){\vector(2,1){.07}}\multiput(4.25,14.25)(.0619800333,.0336938436){601}{\line(1,0){.0619800333}}
\put(41.375,51.875){\vector(0,-1){.07}}\multiput(41.5,69)(-.03125,-4.28125){8}{\line(0,-1){4.28125}}
\put(18.5,48){\makebox(0,0)[cc]{$y_1$}}
\put(64.75,46.75){\makebox(0,0)[cc]{$\bar{y}_1$}}
\put(39.75,07.75){\makebox(0,0)[cc]{$y_0$}}
\put(51.75,22.25){\makebox(0,0)[cc]{$x_0$}}
\put(27,22.5){\makebox(0,0)[cc]{$x_2$}}
\put(36.25,43.25){\makebox(0,0)[cc]{$x_1$}}
\end{picture}
\end{center}

So we can think about the Cayley graph of $F$ in symmetric set of generators $\{y_1,y_1^{-1}\}$. Now we use~\cite[Theorem 1]{Gu21a} where we proved that the density of this Cayley graph is strictly greater than $3$. More precisely, we constructed a finite subgraph $Y$ with density (the average vertex degree) greater than $3$. We do not need to recall all the way of constructing this graph. Just notice that we were based on a family of Brown -- Belk sets $BB(n,k)$ from~\cite{BB05}. They consist of all marked forests with $n$ leaves, where the hight of all their trees does not exceed $k$. We deleted some isolated vertices of Cayley subgraphs for the symmetric generating set $\{x_1,x_1^{-1}\}$. Now these generators are replaced by $\{y_1,y_1^{-1}\}$. This change will be also applied for the construction from~\cite{Gu21a} we are discussing.

The only important property we need to extract now from the proof of the above theorem is the following. The automation $Y$ accepts $y_1$ if and only if the corresponding marked tree of the marked binary forest is nonempty. The same holds for the generator $\bar{y}_1$. Let $N=|Y|$ be the number of vertices in the subgraph $Y$. By $M$ we denote the number of vertices of $Y$ that accept both $y_1$ and $\bar{y}_1$. In this case the number of directed edges of $Y$ labelled by $y_1$ or $\bar{y}_1$ will be equal to $2M$. The number of inverse directed edges will be the same. So the total number of directed edges in the subgraph $Y$ will be equal to $4M$. The density of $Y$ in generators $y_1$, $\bar{y}_1$ is the number of directed edges divided by the number of vertices. So it is $\delta=4M/N$.

According to the construction from~\cite[Theorem 1]{Gu21a}, we have $\delta > 3$. This means that $\frac{N}M < \frac43$. Now let us construct a bipartite graph. For any triangle in $Y$ sided by $y_1$, $\bar{y}_1$, $y_0$ as in the picture, we choose a vertex of the Cayley graph of $F$ in the generating set $A$ inside this triangle. They form a set $S$, where $|S|=M$. If $g\in S$ is one of these vertices, then the vertices of the triangle will be $x_0g$, $x_1g$, $x_2g$ according to the definition of left Cayley graphs. 

Notice that all isolated vertices of $Y$ can be deleted inreasing its density. So we may assume that each vertex of $Y$ is an endpoint of an edge labelled by $y_1$ or $\bar{y}_1$. So it is a vertex of one of the triangles. This means that $AS=Y$. This gives a bipartite graph with independent sets $S$ and $Y$. So we have $|AS|=|Y|=N=\frac{4M}{\delta} < \frac43M=\frac43|S|$. This completes the proof.
\vspace{1ex}

From here we get an immediate

\begin{cy}
\label{qq3}
Let $R=K[F]$ be a group ring of $F$ over a field $K$. For any $4$ linear combinations of elements $x_0$, $x_1$, $x_2$ with coefficients in $K$, the system of $3$ equations with $4$ unknowns
$$
(\alpha_0x_0+\beta_0x_1+\gamma_0x_2)u_0=(\alpha_1x_0+\beta_1x_1+\gamma_1x_2)u_1=(\alpha_2x_0+\beta_2x_1+\gamma_2x_2)u_2=(\alpha_3x_0+\beta_3x_1+\gamma_3x_2)u_3
$$
has a non-zero solution in $R$.
\end{cy}

The proof goes directly using cardinality reasons. We take the set $S$ from Theorem~\ref{43} and express the unkowns $u_0$, \dots, $u_3$ as linear combinations of elements in $S$ with free coefficients from the field. The number of these coefficients is $4|S|$; the number of ordinary linear equations we get from here making each coefficient on an element in $AS$ zero, is $3|AS|$. The inequality $4|S| > 3|AS|$ gives us a non-zero solution of the system of ordinary linear equations.
\vspace{1ex}

Notice that the size of the sets used in the proof of Corollary~\ref{qq3} is very huge. It looks like computer search cannot help to find such a solution explicitly. Let us add that Problem ${\cal Q}_4$ looks like a candidate to a negative answer.

\section{Equations $(1-x_0)u=bv$}
\label{eqx0b}

Throughout the rest of the paper, let $\phi\colon F\to F$ be a shift endomorphism that takes $x_i$ to $x_{i+1}$ for any $i\ge0$. This can be extended to the mapping from the group ring $R=K[F]$ into itself, where $K$ is a field.
\vspace{1ex}

Recall that in~\cite{Gu21b} we reduced the equation problem for $K[F]$ to the case of the monoid ring $K[M]$. So we can work with equations for this situation that looks simpler. We also reduced the equation problem to the case of homogeneous polynomials. However, now we are going to work in a more general situation. 

The main result of this Section is the following

\begin{thm}
\label{x0b}
Let $R=K[F]$ be a group ring of $F$ over a field $K$. Then for any element $b\in R$, the equation $(1-x_0)u=bv$ has a non-zero solution in $R$.
\end{thm}

Notice that this fact is quite general since $b$ can be any element from the group ring. The proof does not go in an obvious way. We start from an auxiliary fact that looks interesting in itself.

\begin{thm}
\label{bphi}
Let $b\in K[M]$ be any element in the monoid ring $R=K[M]$. Then the set $B=\{b,\phi(b),\phi^2(b),\cdots\}$ is not a free basis of the right $R$-module it generates.
\end{thm}

In other words, there exists some integer $m$ such that the equation $bu_0+\phi(b)u_1+\cdots+\phi^m(b)u_m=0$ has a non-zero solution in $R$.

\prf We can assume that $b$ is a sum of homogeneous components: $b=h_1+\cdots+h_k$, where $\deg h_1 < \cdots < \deg h_k$. We proceed by induction on $k$.

Let $k=1$, that is, $b$ is homogeneous. Then $b=f(x_i,\ldots,x_{i+s})$ is a homogeneous polynomial of degree $d$ in some variables. Clearly, $\phi^m(b)=f(x_{i+m},\ldots,x_{i+s+m})$. Assume that $m > s$. Let $w,v\in M$ be any monomials such that each letter $x_j$ occurring in $w$ has a greater subscript than any letter $x_i$ occurring in $v$. Then, according to the defining relations of $M$, the product $wv$ equals $v\phi^d(w)$, where $d$ is the length (degree) of $v$. The same effect holds if we replace $v$ be a homogeneous polynomial of degree $d$, and replace $w$ by any polynomial. We only claim that any letter $x_j$ involved in $w$ has a greater subscript than any letter $x_i$ involved in $v$. In this situation we also have $wv=v\phi^d(w)$.

Applying this to our situation, we get $\phi^m(b)b=f(x_{i+m},\ldots,x_{i+s+m})f(x_i,\ldots,x_{i+s})=f(x_i,\ldots,x_{i+s})f(x_{i+m+d},\ldots,x_{i+s+m+d})=b\phi^{m+d}(b)$. This shows that $bR$ and $\phi^m(b)R$ have non-zero intersection, which is enough for our case.

Now let $k > 1$. Here we also regard $b$ as a polynomial in $x_i$, \dots, $x_{i+s}$. Let $m > s$ and $d=\deg h_k$. Then 
$$
\phi^m(b)\cdot h_k=h_k\phi^{m+d}(b)=h_k\phi^{m+d}(h_1)+\cdots+h_k\phi^{m+d}(h_{k-1})+h_k\phi^{m+d}(h_k);
$$
$$
b\cdot\phi^{m+d}(h_k)=h_1\phi^{m+d}(h_k)+\cdots+h_{k-1}\phi^{m+d}(h_k)+h_k\phi^{m+d}(h_k).
$$
Then $c\leftrightharpoons\phi^m(b)h_k-b\phi^{m+d}(h_k)=h_1'+\cdots+h_{k-1}'$, where $h_j'=h_k\phi^{m+d}(h_j)-h_j\phi^{m+d}(h_k)$ are homogeneous polynomials of degree $\deg h_j+d$ ($1\le j\le k-1$). So $c$ is the element from the right ideal generated by $B$, being a sum of less than $k$ homogeneous components. By the inductive assumption, there exists a nontrivial right module relation between $c$, $\phi(c)$, $\phi^2(c)$, \dots\,. This leads to a nontrivial module relation between elements of $B$.
\vspace{1ex}

{\bf Proof of Theorem~\ref{x0b}.} Due to Lemma~\ref{gig}, multiplying by elements from $F$, we reduce the problem to the case $b\in K[M]$ finding solutions in the monoid ring $K[M]$. Let $b=b_0+x_0b_1+\cdots+x_0^mb_m$, where $b_0,b_1,\ldots,b_m$ do not contain $x_0$. That is, they belong to the monoid ring $K[M_1]$, where $M_1$ is the submonoid of $M$ generated by $x_1$, $x_2$, \dots\,. Let $b'=b_0+b_1+\cdots+b_m$. Hence $b'-b=(1-x_0)b_1+\cdots+(1-x_0^m)b_m$ belongs to $(1-x_0)R$, that is, $b'=b+(1-x_0)c$ for some $c\in K[M]$.

It suffices to prove that equation $(1-x_0)u=b'v$ has a nonzero solution. Indeed, this is equivalent to $(1-x_0)u=bv+(1-x_0)cv$, so $(1-x_0)(u-cv)=bv$ gives a non-zero solution of the original equation. Notice that $b'\in K[M_1]$ does not depend on $x_0$. We are going to find $v=v_0+x_0v_1+\cdots+x_0^kv_k$, where $v_0,v_1,\ldots,v_k\in K[M_1]$ and the element $b'v$ is left divisible by $1-x_0$. We have $b'v=b'v_0+x_0\phi(b')v_1+x_0^2\phi^2(b')v_2+\cdots+x_0^k\phi^k(b')v_k$. 

Let us apply Theorem~\ref{bphi} to the set $B'=\{b',\phi(b'),\phi^2(b'),\ldots\}$. This gives us a non-zero solution to the equation $b'v_0+\phi(b')v_1+\cdots+\phi^k(b')v_k=0$ for some $k$. Therefore, $b'v=(x_0-1)\phi(b')v_1+\cdots+(x_0^k-1)\phi^k(b')v_k$ is left divisible by $1-x_0$ in the ring, that is, $b'v=(1-x_0)u$ for some $u\in R$. Obviously, $v\ne0$ since at least one of the $v_i$s is non-zero. So we have a nonzero solution of the equation $(1-x_0)u=b'v$, as desired.

The proof is complete.
\vspace{1ex}

We do not know the answer to the same question on equations from Theorem~\ref{x0b} if we replace $1-x_0$ by $1-x_1$. So we ask the following

\begin{que}
\label{x1b}
Let $R=K[F]$ be a group ring of $F$ over a field $K$. Is it true that for any element $b\in R$, the equation $(1-x_1)u=bv$ has a non-zero solution in $R$?
\end{que}

Notice that applying the involutive automorphism $x_0\mapsto x_0^{-1}$, $x_1\mapsto\bar{x}_1=x_1x_0^{-1}$, one can reduce the above question to the equation of the form $(x_0-x_1)u=bv$. Without loss of generality, $b\in K[M]$. One can expect a negative answer for a homogeneous polynomial $b$ of sufficiently high degree.
\vspace{1ex}

Now we are going to deduce a corollary from Theorem~\ref{x0b}. Let us find a non-zero solution for $(1-x_0)u=(1-x_1)v$. Applying $\phi$, we get the element $b=(1-x_1)\phi(u)=(1-x_2)\phi(v)$. Changing notation and solving equation $(1-x_0)u=bv$ according to Theorem~\ref{x0b}, we get a nonzero element of the form $(1-x_0)u=(1-x_1)v=(1-x_2)w$ and so on. Applying $\phi$ and the same Theorem, we get by induction the following

\begin{cy}
\label{x0xm}
For any $m\ge1$ there exists a non-zero solution to the system of equations
\be{x0x1xm}
(1-x_0)u_0=(1-x_1)u_1=\cdots=(1-x_m)u_m
\ee
in the group ring of $F$.
\end{cy}

Of course, one can present such a solution in the integer ring $\mathbb Z[M]$. The support $G_m$ of the element in~(\ref{x0x1xm}) promises to have interesting properties. We can expect the following effect. An integer linear combination of elements of $G_m$ is left divisible by any $1-x_i$, where $0\le i\le m$. This usually happens whenever for any $i$ the elements of $G_m$ can be paired in such a way that two elements forming a pair differ by a power of $x_i$. The case $m=1$ is just any relation in $F$; cf. the beginning of the next Section. So one can regard the elements from~(\ref{x0x1xm}) as multi-dimensional analogs of relations between generators.

We did not find an explicit form of~(\ref{x0x1xm}) even for the case $m=2$. It is interesting to imagine how this solution will look like. It does not seem to have a very big size so a computer search can help to find minimal solutions of the system for small values of $m$.

\section{Solutions of the equation $(1-x_0)u=(1-x_1)v$}
\label{eqx0x1}

In this Section we present an explicit description for the set of all solutions of the given equation in the group ring of $F$. Of course, it is enough to describe this set for the monoid ring $K[M]$. We start from an easy observation that shows a close  connection between relations in $F$ and solutions to the equation in the title of this Section.

In the picture below we take one of the defining relations in $F$, namely, $x_1^{x_0^2}=x_1^{x_0x_1}$. We regard it as a part of the left Cayley graph. By a right multiplication, each relation can be placed inside $M$ according to Lemma~\ref{gig}. 

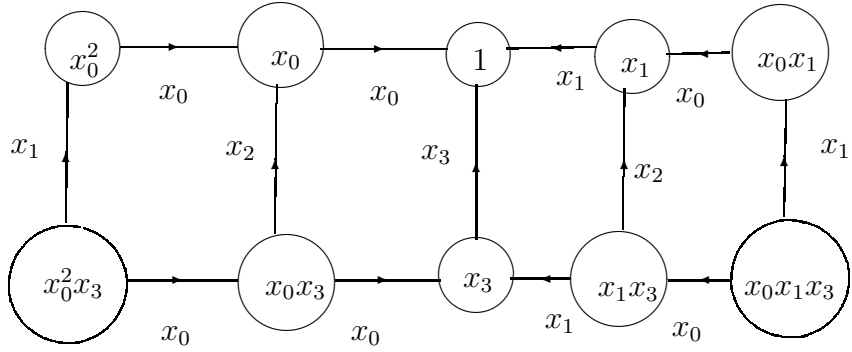
\begin{center}
	\begin{picture}(112.766,46.119)(0,0)
	\put(12.25,38.5){\makebox(0,0)[cc]{$x_0^2$}}
	\put(12,40){\circle{10.198}}
	\put(38.75,38.5){\makebox(0,0)[cc]{$x_0$}}
	\put(38,40.5){\circle{10.548}}
	\put(64.25,38.5){\makebox(0,0)[cc]{$1$}}
	\put(64,39.25){\circle{8.062}}
	\put(84.75,37.5){\makebox(0,0)[cc]{$x_1$}}
	\put(84.25,39){\circle{9.179}}
	\put(104.75,38.25){\makebox(0,0)[cc]{$x_0x_1$}}
	\put(103.75,39.5){\circle{13.238}}
	\put(10.75,8.75){\makebox(0,0)[cc]{$x_0^2x_3$}}
	\put(17.766,8.75){\line(0,1){.4373}}
	\put(17.754,9.187){\line(0,1){.4359}}
	\put(17.717,9.623){\line(0,1){.4331}}
	\put(17.655,10.056){\line(0,1){.429}}
	\put(17.57,10.485){\line(0,1){.4235}}
	\multiput(17.46,10.909)(-.03335,.10415){4}{\line(0,1){.10415}}
	\multiput(17.327,11.325)(-.031326,.081689){5}{\line(0,1){.081689}}
	\multiput(17.17,11.734)(-.029896,.066496){6}{\line(0,1){.066496}}
	\multiput(16.991,12.133)(-.033593,.064708){6}{\line(0,1){.064708}}
	\multiput(16.789,12.521)(-.031871,.053754){7}{\line(0,1){.053754}}
	\multiput(16.566,12.897)(-.030491,.04539){8}{\line(0,1){.04539}}
	\multiput(16.322,13.26)(-.032999,.043601){8}{\line(0,1){.043601}}
	\multiput(16.058,13.609)(-.031468,.037044){9}{\line(0,1){.037044}}
	\multiput(15.775,13.943)(-.033504,.035213){9}{\line(0,1){.035213}}
	\multiput(15.473,14.26)(-.035433,.033271){9}{\line(-1,0){.035433}}
	\multiput(15.154,14.559)(-.03725,.031223){9}{\line(-1,0){.03725}}
	\multiput(14.819,14.84)(-.043818,.032711){8}{\line(-1,0){.043818}}
	\multiput(14.469,15.102)(-.04559,.030192){8}{\line(-1,0){.04559}}
	\multiput(14.104,15.343)(-.053963,.031516){7}{\line(-1,0){.053963}}
	\multiput(13.726,15.564)(-.064928,.033166){6}{\line(-1,0){.064928}}
	\multiput(13.337,15.763)(-.066692,.029457){6}{\line(-1,0){.066692}}
	\multiput(12.936,15.94)(-.081894,.030787){5}{\line(-1,0){.081894}}
	\multiput(12.527,16.093)(-.10437,.03266){4}{\line(-1,0){.10437}}
	\put(12.11,196.44){\line(-1,0){.4242}}
	\put(11.685,16.331){\line(-1,0){.4295}}
	\put(11.256,16.414){\line(-1,0){.4335}}
	\put(10.822,16.472){\line(-1,0){.4361}}
	\put(10.386,16.507){\line(-1,0){.8745}}
	\put(9.512,16.501){\line(-1,0){.4356}}
	\put(9.076,16.461){\line(-1,0){.4327}}
	\put(8.643,16.397){\line(-1,0){.4284}}
	\put(8.215,16.308){\line(-1,0){.4227}}
	\multiput(7.792,16.196)(-.083146,-.027225){5}{\line(-1,0){.083146}}
	\multiput(7.376,16.06)(-.081481,-.031863){5}{\line(-1,0){.081481}}
	\multiput(6.969,15.9)(-.066298,-.030334){6}{\line(-1,0){.066298}}
	\multiput(6.571,15.718)(-.055273,-.029159){7}{\line(-1,0){.055273}}
	\multiput(6.184,15.514)(-.053543,-.032225){7}{\line(-1,0){.053543}}
	\multiput(5.81,15.289)(-.045188,-.03079){8}{\line(-1,0){.045188}}
	\multiput(5.448,15.042)(-.043383,-.033285){8}{\line(-1,0){.043383}}
	\multiput(5.101,14.776)(-.036836,-.031711){9}{\line(-1,0){.036836}}
	\multiput(4.769,14.491)(-.034992,-.033735){9}{\line(-1,0){.034992}}
	\multiput(4.455,14.187)(-.033037,-.035652){9}{\line(0,-1){.035652}}
	\multiput(4.157,13.866)(-.030977,-.037455){9}{\line(0,-1){.037455}}
	\multiput(3.878,13.529)(-.032421,-.044032){8}{\line(0,-1){.044032}}
	\multiput(3.619,13.177)(-.029891,-.045788){8}{\line(0,-1){.045788}}
	\multiput(3.38,12.81)(-.03116,-.05417){7}{\line(0,-1){.05417}}
	\multiput(3.162,12.431)(-.032737,-.065145){6}{\line(0,-1){.065145}}
	\multiput(2.965,12.04)(-.029017,-.066885){6}{\line(0,-1){.066885}}
	\multiput(2.791,11.639)(-.030247,-.082095){5}{\line(0,-1){.082095}}
	\multiput(2.64,11.229)(-.03197,-.10458){4}{\line(0,-1){.10458}}
	\put(2.512,10.81){\line(0,-1){.4249}}
	\put(2.408,10.385){\line(0,-1){.4301}}
	\put(2.328,9.955){\line(0,-1){.4339}}
	\put(2.272,9.521){\line(0,-1){.8737}}
	\put(2.235,8.648){\line(0,-1){.4371}}
	\put(2.253,8.211){\line(0,-1){.4354}}
	\put(2.295,7.775){\line(0,-1){.4323}}
	\put(2.362,7.343){\line(0,-1){.4278}}
	\put(2.454,6.915){\line(0,-1){.422}}
	\multiput(2.569,6.493)(.027772,-.082965){5}{\line(0,-1){.082965}}
	\multiput(2.708,6.078)(.032399,-.081269){5}{\line(0,-1){.081269}}
	\multiput(2.87,5.672)(.03077,-.066097){6}{\line(0,-1){.066097}}
	\multiput(3.055,5.275)(.029522,-.055079){7}{\line(0,-1){.055079}}
	\multiput(3.261,4.89)(.032577,-.05333){7}{\line(0,-1){.05333}}
	\multiput(3.489,4.517)(.031087,-.044985){8}{\line(0,-1){.044985}}
	\multiput(3.738,4.157)(.03357,-.043163){8}{\line(0,-1){.043163}}
	\multiput(4.006,3.811)(.031953,-.036626){9}{\line(0,-1){.036626}}
	\multiput(4.294,3.482)(.030568,-.031292){10}{\line(0,-1){.031292}}
	\multiput(4.6,3.169)(.035869,-.032801){9}{\line(1,0){.035869}}
	\multiput(4.923,2.874)(.037659,-.03073){9}{\line(1,0){.037659}}
	\multiput(5.261,2.597)(.044245,-.032131){8}{\line(1,0){.044245}}
	\multiput(5.615,2.34)(.045984,-.029588){8}{\line(1,0){.045984}}
	\multiput(5.983,2.103)(.054374,-.030803){7}{\line(1,0){.054374}}
	\multiput(6.364,1.888)(.065359,-.032308){6}{\line(1,0){.065359}}
	\multiput(6.756,1.694)(.067074,-.028576){6}{\line(1,0){.067074}}
	\multiput(7.159,1.522)(.082292,-.029705){5}{\line(1,0){.082292}}
	\multiput(7.57,1.374)(.10479,-.03128){4}{\line(1,0){.10479}}
	\put(7.989,1.249){\line(1,0){.4256}}
	\put(8.415,1.147){\line(1,0){.4306}}
	\put(8.845,1.07){\line(1,0){.4342}}
	\put(9.28,1.017){\line(1,0){.4365}}
	\put(9.716,0.989){\line(1,0){.4374}}
	\put(10.153,0.985){\line(1,0){.4369}}
	\put(10.59,1.006){\line(1,0){.4351}}
	\put(11.025,1.052){\line(1,0){.4318}}
	\put(11.457,1.122){\line(1,0){.4272}}
	\put(11.885,1.216){\line(1,0){.4212}}
	\multiput(12.306,1.334)(.08278,.028318){5}{\line(1,0){.08278}}
	\multiput(12.72,1.476)(.081054,.032934){5}{\line(1,0){.081054}}
	\multiput(13.125,1.64)(.065893,.031205){6}{\line(1,0){.065893}}
	\multiput(13.52,1.828)(.054884,.029884){7}{\line(1,0){.054884}}
	\multiput(13.904,2.037)(.053114,.032927){7}{\line(1,0){.053114}}
	\multiput(14.276,2.267)(.044779,.031382){8}{\line(1,0){.044779}}
	\multiput(14.634,2.518)(.03817,.030092){9}{\line(1,0){.03817}}
	\multiput(14.978,2.789)(.036415,.032194){9}{\line(1,0){.036415}}
	\multiput(15.306,3.079)(.03109,.030774){10}{\line(1,0){.03109}}
	\multiput(15.617,3.387)(.032564,.036084){9}{\line(0,1){.036084}}
	\multiput(15.91,3.711)(.030481,.03786){9}{\line(0,1){.03786}}
	\multiput(16.184,4.052)(.031838,.044456){8}{\line(0,1){.044456}}
	\multiput(16.439,4.408)(.033468,.052775){7}{\line(0,1){.052775}}
	\multiput(16.673,4.777)(.030444,.054576){7}{\line(0,1){.054576}}
	\multiput(16.886,5.159)(.031876,.06557){6}{\line(0,1){.06557}}
	\multiput(17.077,5.553)(.028134,.067261){6}{\line(0,1){.067261}}
	\multiput(17.246,5.956)(.029162,.082486){5}{\line(0,1){.082486}}
	\multiput(17.392,6.369)(.03059,.105){4}{\line(0,1){.105}}
	\put(17.514,6.789){\line(0,1){.4262}}
	\put(17.613,7.215){\line(0,1){.4311}}
	\put(17.687,7.646){\line(0,1){.4346}}
	\put(17.737,8.08){\line(0,1){.6695}}
	\put(40,8){\makebox(0,0)[cc]{$x_0x_3$}}
	\put(38.75,9){\circle{12.816}}
	\put(64,9){\makebox(0,0)[cc]{$x_3$}}
	\put(63.75,10){\circle{9.394}}
	\put(83.75,7.75){\makebox(0,0)[cc]{$x_1x_3$}}
	\put(82.5,9.5){\circle{12.855}}
	\put(105,8){\makebox(0,0)[cc]{$x_0x_1x_3$}}
	\put(112.766,9.5){\line(0,1){.4373}}
	\put(112.754,9.937){\line(0,1){.4359}}
	\put(112.717,10.373){\line(0,1){.4331}}
	\put(112.655,10.806){\line(0,1){.429}}
	\put(112.57,11.235){\line(0,1){.4235}}
	\multiput(112.46,11.659)(-.03335,.10415){4}{\line(0,1){.10415}}
	\multiput(112.327,12.075)(-.031326,.081689){5}{\line(0,1){.081689}}
	\multiput(112.17,12.484)(-.029896,.066496){6}{\line(0,1){.066496}}
	\multiput(111.991,12.883)(-.033593,.064708){6}{\line(0,1){.064708}}
	\multiput(111.789,13.271)(-.031871,.053754){7}{\line(0,1){.053754}}
	\multiput(111.566,13.647)(-.030491,.04539){8}{\line(0,1){.04539}}
	\multiput(111.322,14.01)(-.032999,.043601){8}{\line(0,1){.043601}}
	\multiput(111.058,14.359)(-.031468,.037044){9}{\line(0,1){.037044}}
	\multiput(110.775,14.693)(-.033504,.035213){9}{\line(0,1){.035213}}
	\multiput(110.473,15.01)(-.035433,.033271){9}{\line(-1,0){.035433}}
	\multiput(110.154,15.309)(-.03725,.031223){9}{\line(-1,0){.03725}}
	\multiput(109.819,15.59)(-.043818,.032711){8}{\line(-1,0){.043818}}
	\multiput(109.469,15.852)(-.04559,.030192){8}{\line(-1,0){.04559}}
	\multiput(109.104,16.093)(-.053963,.031516){7}{\line(-1,0){.053963}}
	\multiput(108.726,16.314)(-.064928,.033166){6}{\line(-1,0){.064928}}
	\multiput(108.337,16.513)(-.066692,.029457){6}{\line(-1,0){.066692}}
	\multiput(107.936,16.69)(-.081894,.030787){5}{\line(-1,0){.081894}}
	\multiput(107.527,16.843)(-.10437,.03266){4}{\line(-1,0){.10437}}
	\put(107.11,16.974){\line(-1,0){.4242}}
	\put(106.685,17.081){\line(-1,0){.4295}}
	\put(106.256,17.164){\line(-1,0){.4335}}
	\put(105.822,17.222){\line(-1,0){.4361}}
	\put(105.386,17.257){\line(-1,0){.8745}}
	\put(104.512,17.251){\line(-1,0){.4356}}
	\put(104.076,17.211){\line(-1,0){.4327}}
	\put(103.643,17.147){\line(-1,0){.4284}}
	\put(103.215,17.058){\line(-1,0){.4227}}
	\multiput(102.792,16.946)(-.083146,-.027225){5}{\line(-1,0){.083146}}
	\multiput(102.376,16.81)(-.081481,-.031863){5}{\line(-1,0){.081481}}
	\multiput(101.969,16.65)(-.066298,-.030334){6}{\line(-1,0){.066298}}
	\multiput(101.571,16.468)(-.055273,-.029159){7}{\line(-1,0){.055273}}
	\multiput(101.184,16.264)(-.053543,-.032225){7}{\line(-1,0){.053543}}
	\multiput(100.81,16.039)(-.045188,-.03079){8}{\line(-1,0){.045188}}
	\multiput(100.448,15.792)(-.043383,-.033285){8}{\line(-1,0){.043383}}
	\multiput(100.101,15.526)(-.036836,-.031711){9}{\line(-1,0){.036836}}
	\multiput(99.769,15.241)(-.034992,-.033735){9}{\line(-1,0){.034992}}
	\multiput(99.455,14.937)(-.033037,-.035652){9}{\line(0,-1){.035652}}
	\multiput(99.157,14.616)(-.030977,-.037455){9}{\line(0,-1){.037455}}
	\multiput(98.878,14.279)(-.032421,-.044032){8}{\line(0,-1){.044032}}
	\multiput(98.619,13.927)(-.029891,-.045788){8}{\line(0,-1){.045788}}
	\multiput(98.38,13.56)(-.03116,-.05417){7}{\line(0,-1){.05417}}
	\multiput(98.162,13.181)(-.032737,-.065145){6}{\line(0,-1){.065145}}
	\multiput(97.965,12.79)(-.029017,-.066885){6}{\line(0,-1){.066885}}
	\multiput(97.791,12.389)(-.030247,-.082095){5}{\line(0,-1){.082095}}
	\multiput(97.64,11.979)(-.03197,-.10458){4}{\line(0,-1){.10458}}
	\put(97.512,11.56){\line(0,-1){.4249}}
	\put(97.408,11.135){\line(0,-1){.4301}}
	\put(97.328,10.705){\line(0,-1){.4339}}
	\put(97.272,10.271){\line(0,-1){.8737}}
	\put(97.235,9.398){\line(0,-1){.4371}}
	\put(97.253,8.961){\line(0,-1){.4354}}
	\put(97.295,8.525){\line(0,-1){.4323}}
	\put(97.362,8.093){\line(0,-1){.4278}}
	\put(97.454,7.665){\line(0,-1){.422}}
	\multiput(97.569,7.243)(.027772,-.082965){5}{\line(0,-1){.082965}}
	\multiput(97.708,6.828)(.032399,-.081269){5}{\line(0,-1){.081269}}
	\multiput(97.87,6.422)(.03077,-.066097){6}{\line(0,-1){.066097}}
	\multiput(98.055,6.025)(.029522,-.055079){7}{\line(0,-1){.055079}}
	\multiput(98.261,5.64)(.032577,-.05333){7}{\line(0,-1){.05333}}
	\multiput(98.489,5.267)(.031087,-.044985){8}{\line(0,-1){.044985}}
	\multiput(98.738,4.907)(.03357,-.043163){8}{\line(0,-1){.043163}}
	\multiput(99.006,4.561)(.031953,-.036626){9}{\line(0,-1){.036626}}
	\multiput(99.294,4.232)(.030568,-.031292){10}{\line(0,-1){.031292}}
	\multiput(99.6,3.919)(.035869,-.032801){9}{\line(1,0){.035869}}
	\multiput(99.923,3.624)(.037659,-.03073){9}{\line(1,0){.037659}}
	\multiput(100.261,3.347)(.044245,-.032131){8}{\line(1,0){.044245}}
	\multiput(100.615,3.09)(.045984,-.029588){8}{\line(1,0){.045984}}
	\multiput(100.983,2.853)(.054374,-.030803){7}{\line(1,0){.054374}}
	\multiput(101.364,2.638)(.065359,-.032308){6}{\line(1,0){.065359}}
	\multiput(101.756,2.444)(.067074,-.028576){6}{\line(1,0){.067074}}
	\multiput(102.159,2.272)(.082292,-.029705){5}{\line(1,0){.082292}}
	\multiput(102.57,2.124)(.10479,-.03128){4}{\line(1,0){.10479}}
	\put(102.989,1.999){\line(1,0){.4256}}
	\put(103.415,1.897){\line(1,0){.4306}}
	\put(103.845,1.82){\line(1,0){.4342}}
	\put(104.28,1.767){\line(1,0){.4365}}
	\put(104.716,1.739){\line(1,0){.4374}}
	\put(105.153,1.735){\line(1,0){.4369}}
	\put(105.59,1.756){\line(1,0){.4351}}
	\put(106.025,1.802){\line(1,0){.4318}}
	\put(106.457,1.872){\line(1,0){.4272}}
	\put(106.885,1.966){\line(1,0){.4212}}
	\multiput(107.306,2.084)(.08278,.028318){5}{\line(1,0){.08278}}
	\multiput(107.72,2.226)(.081054,.032934){5}{\line(1,0){.081054}}
	\multiput(108.125,2.39)(.065893,.031205){6}{\line(1,0){.065893}}
	\multiput(108.52,2.578)(.054884,.029884){7}{\line(1,0){.054884}}
	\multiput(108.904,2.787)(.053114,.032927){7}{\line(1,0){.053114}}
	\multiput(109.276,3.017)(.044779,.031382){8}{\line(1,0){.044779}}
	\multiput(109.634,3.268)(.03817,.030092){9}{\line(1,0){.03817}}
	\multiput(109.978,3.539)(.036415,.032194){9}{\line(1,0){.036415}}
	\multiput(110.306,3.829)(.03109,.030774){10}{\line(1,0){.03109}}
	\multiput(110.617,4.137)(.032564,.036084){9}{\line(0,1){.036084}}
	\multiput(110.91,4.461)(.030481,.03786){9}{\line(0,1){.03786}}
	\multiput(111.184,4.802)(.031838,.044456){8}{\line(0,1){.044456}}
	\multiput(111.439,5.158)(.033468,.052775){7}{\line(0,1){.052775}}
	\multiput(111.673,5.527)(.030444,.054576){7}{\line(0,1){.054576}}
	\multiput(111.886,5.909)(.031876,.06557){6}{\line(0,1){.06557}}
	\multiput(112.077,6.303)(.028134,.067261){6}{\line(0,1){.067261}}
	\multiput(112.246,6.706)(.029162,.082486){5}{\line(0,1){.082486}}
	\multiput(112.392,7.119)(.03059,.105){4}{\line(0,1){.105}}
	\put(112.514,7.539){\line(0,1){.4262}}
	\put(112.613,7.965){\line(0,1){.4311}}
	\put(112.687,8.396){\line(0,1){.4346}}
	\put(112.737,8.83){\line(0,1){.6695}}
	\put(24.625,40.25){\vector(1,0){.07}}\put(17,40.25){\line(1,0){15.25}}
	\put(51.625,40){\vector(1,0){.07}}\put(43.25,40){\line(1,0){16.75}}
	\put(25.125,9.25){\vector(1,0){.07}}\put(18,9.25){\line(1,0){14.25}}
	\put(52.125,9.25){\vector(1,0){.07}}\put(45.25,9.25){\line(1,0){13.75}}
	\put(9.75,26.25){\vector(0,1){.07}}\multiput(9.75,17)(.03125,2.3125){8}{\line(0,1){2.3125}}
	\put(37.275,25.375){\vector(0,1){.07}}\multiput(37.25,15.5)(.03125,2.46875){8}{\line(0,1){2.46875}}
	\put(63.75,25){\vector(0,1){.07}}\put(63.75,14.75){\line(0,1){20.5}}
	\put(83.16,25.375){\vector(0,1){.07}}\multiput(83,16)(.03125,2.34375){8}{\line(0,1){2.34375}}
	\put(104.25,25.75){\vector(0,1){.07}}\multiput(104.25,17.75)(.03125,1.9375){8}{\line(0,1){1.9375}}
	\put(72.125,9.5){\vector(-1,0){.07}}\put(76,9.5){\line(-1,0){7.75}}
	\put(93,9.25){\vector(-1,0){.07}}\put(97,9.25){\line(-1,0){8}}
	\put(73.625,40.125){\vector(-1,0){.07}}\multiput(79.5,40.25)(-1.46875,-.03125){8}{\line(-1,0){1.46875}}
	\put(92.875,39.375){\vector(-1,0){.07}}\multiput(97,39.5)(-1.03125,-.03125){8}{\line(-1,0){1.03125}}
	\put(24.25,2){\makebox(0,0)[cc]{$x_0$}}
	\put(49.25,2){\makebox(0,0)[cc]{$x_0$}}
	\put(91.5,2.5){\makebox(0,0)[cc]{$x_0$}}
	\put(92,33.5){\makebox(0,0)[cc]{$x_0$}}
	\put(51.75,33.75){\makebox(0,0)[cc]{$x_0$}}
	\put(24,34){\makebox(0,0)[cc]{$x_0$}}
	\put(4.5,26.75){\makebox(0,0)[cc]{$x_1$}}
	\put(111,26.75){\makebox(0,0)[cc]{$x_1$}}
	\put(76.25,35.75){\makebox(0,0)[cc]{$x_1$}}
	\put(74.75,3.5){\makebox(0,0)[cc]{$x_1$}}
	\put(32.75,26.5){\makebox(0,0)[cc]{$x_2$}}
	\put(58.5,26){\makebox(0,0)[cc]{$x_3$}}
	\put(86.5,23.5){\makebox(0,0)[cc]{$x_2$}}
	\end{picture}
\end{center}

Going from the vertex $1$ counterclockwise and taking into account that $x_0^2x_3=x_1x_0^2$, $x_0x_1x_3=x_1x_0x_1$ in $F$, we write $0=1-x_0+x_0-x_0^2+x_0^2-x_1x_0^2+x_0^2x_3-x_0x_3+x_0x_3-x_3+x_3-x_1x_3+x_1x_3-x_0x_1x_3+x_1x_0x_1-x_0x_1+x_0x_1-x_1+x_1-1$. This means that $0=(1-x_0)(1+x_0-x_0x_3-x_3+x_1x_3-x_1)+(1-x_1)(x_0^2+x_3-x_0x_1-1)$. It follows that $(1-x_0)u_0=(1-x_1)v_0$, where
\be{u0v0}
u_0=(1+x_0-x_1)(1-x_3),\ \ \ v_0=1-x_3-x_0^2+x_0x_1.
\ee
We call the pair $(u_0,v_0)$ a {\em basic solution} to the equation $(1-x_0)u=(1-x_1)v$.

Now we are going to extract all solutions of this equation from the basic one. Clearly, we have solutions of the form $u=u_0w$, $v=v_0w$ for any $w\in K[M]$. 
\vspace{1ex}

For any $i\ge0$ by $M_i$ we denote the submonoid of $M$ generated by $x_i$, $x_{i+1}$, \dots\,. In order to analyze all solutions, we prove the statement similar to~\cite[Lemma 4]{Gu21b}. The idea of it is to divide any element $v\in K[M]$ by $v_0$ with a remainder that does not involve high powers of $x_0$. However, there is one additional condition we need to take into account. In~\cite{Gu21b} we worked with homogeneous polynomials. Now this is not the case so we need to add one restriction on the degrees.

\begin{lm}
\label{v0R}
For any $v\in K[M]$, $v\ne0$ there exist $w_1\in K[M]$, $w_2,w_3\in K[M_1]$ such that $v=v_0w_1+x_0w_2+w_3$. In addition, we can assume that $\deg(x_0w_2),\deg w_3\le\deg v$.
\end{lm}

\begin{prf}
We are going to show that $v=x_0w_2+w_3$ modulo the principal right ideal $v_0R$, where $R=K[M]$. All equalities below will be done modulo this ideal. Each of them can be multiplied on the right.
	
We know that $x_0^2=1-x_3+x_0x_1$. Let us prove by induction on $k\ge2$ that $x_0^k=x_0\xi_k+\eta_k$ for some $\xi_k,\eta_k\in K[M_1]$. In addition, we claim that $\deg\xi_k,\deg\eta_k\le k-1$. We know this for $k=2$, where $\xi_2=x_1$, $\eta_2=1-x_3$. Let this equality and inequalities hold for some $k\ge2$. Then $x_0^{k+1}=(x_0\xi_k+\eta_k)x_0=x_0^2\phi(\xi_k)+x_0\phi(\eta_k)$. Therefore, $x_0^{k+1}=(x_0\xi_2+\eta_2)\phi(\xi_k)+x_0\phi(\eta_k)=x_0\xi_{k+1}+\eta_{k+1}$, where $\xi_{k+1}=\xi_2\phi(\xi_k)+\phi(\eta_k)$, $\eta_{k+1}=\eta_2\phi(\xi_k)$. Clearly, $\phi$ does not change the degree so $\deg\xi_{k+1},\deg\eta_{k+1}\le k$. 
	
Now we can decompose $v$ by powers of $x_0$, that is, $v=\zeta_0+x_0\zeta_1+x_0^2\zeta_2+\cdots+x_0^d\zeta_d$ for some $d$, where $\zeta_0,\ldots,\zeta_d\in K[M_1]$. Replacing here all $x_0^k$ by $x_0\xi_k+\eta_k$ for $2\le k\le d$, we get the desired equality of the form $v=x_0w_2+w_3$ modulo the ideal $v_0R$, where $w_2,w_3\in K[M_1]$. Hence there exists $w_1\in R$ such that $v=v_0w_1+x_0w_2+w_3$. No monomials of degree $>k$ appear in the process of replacements $x_0^k\to x_0\xi_k+\eta_k$. Therefore we never get monomials of degree higher than $\deg v$ in polynomials $x_0w_2$, $w_3$.
	
The proof is complete.
\end{prf}
\vspace{1ex}

Let $(u,v)$ be an arbitrary non-zero solution of the equation $(1-x_0)u=(1-x_1)v$. The basic idea will be to proceed by induction on $\deg v$. Applying Lemma~\ref{v0R}, we set $v=v_0w_1+x_0w_2+w_3$, where $w_1\in K[M]$, $w_2,w_3\in K[M_1]$; $\deg(x_0w_2),\deg w_3\le\deg v$. We know that $(1-x_0)u_0=(1-x_1)v_0$ for our basic solution so $(1-x_0)u_0w_1=(1-x_1)v_0w_1$. Subtracting this equation from $(1-x_0)u=(1-x_1)(v_0w_1+x_0w_2+w_3)$, we get $(1-x_0)(u-u_0w_1)=(1-x_1)(w_3+x_0w_2)$. Let $w=u-u_0w_1$ be a reduced form of the polynomial. We claim that $w\in K[M_1]$. Otherwise $w$ has a monomial containing $x_0$; we take such a monomial of the form $x_0^k\mu$ for the highest $k\ge1$, where $\mu\in K[M_1]$. In the product $(1-x_0)w$ we will get a monomial $x_0^{k+1}\mu$ with nonzero coefficient. On the right-hand side we do not have monomials containg $x_0^2$ since $w_2,w_3$ do not contain $x_0$.
\vspace{1ex}

Now we have $(1-x_0)w=(1-x_1)(w_3+x_0w_2)$. In the left-hand side we have $w-x_0w$ with $w\in K[M_1]$. In the right-hand side we have $(1-x_1)w_3+x_0(1-x_2)w_2$, where $(1-x_1)w_3,(1-x_2)w_2\in K[M_1]$. Comparing both expressions and their parts that belong to $K[M_1]$ and $x_0K[M_1]$, we obtain that $w=(1-x_1)w_3$ and $-w=(1-x_2)w_2$. This gives us an equation $(1-x_1)w_3=(1-x_2)(-w_2)$. Both sides are from $K[M_1]$ so we can apply $\phi^{-1}$ to them getting the equation $(1-x_0)\phi^{-1}(w_3)=(1-x_1)(-\phi^{-1}(w_2))$. 

We know from Lemma~\ref{v0R} that $\deg(x_0w_2)\le\deg v$. Therefore, we have $\deg w_2 < \deg v$. Since $\phi^{-1}$ has no influence on degrees, we potentially know the solutions of the new equation. Namely, if $(u',v')$ is such a solution, then we can express $u$, $v$ from that. Here $u'=\phi^{-1}(w_3)$, $v'=-\phi^{-1}(w_2)$. Hence $w_2=-\phi(v')$, $w_3=\phi(u')$, $w=(1-x_1)w_3=(1-x_1)\phi(u')$, $u=u_0w_1+w=u_0w_1+(1-x_1)\phi(u')$, $v=v_0w_1-x_0\phi(v')+\phi(u')$. 
\vspace{1ex}

Let us summarize this information in the following

\begin{lm}
\label{summ}
Let $u_0=(1+x_0-x_1)(1-x_3)$, $v_0=1-x_3-x_0^2+x_0x_1$ be the basic solution to the equation $(1-x_0)u=(1-x_1)v$. Then for any solution $(u',v')$ of the same equation and for any $w\in K[M]$, the pair 
\be{uuvv}
u=u_0w+(1-x_1)\phi(u'),\ \ \  v=v_0w+\phi(u')-x_0\phi(v')
\ee 
will be a solution to $(1-x_0)u=(1-x_1)v$. 

Moreover, for any non-zero solution $(u,v)$ there exists a solution $(u',v')$ of the same equation with $\deg v' < \deg v$ such that equalities~(\ref{uuvv}) hold for some $w\in K[M]$.
\end{lm}

Now we can present an explicit description for the set of all solutions. Notice that in~(\ref{uuvv}) the equality for $u$ depends of $u'$ only, while $v$ is expressed in terms of both $u'$, $v'$. If we know all possible values of $u$ for the equation $(1-x_0)u=(1-x_1)v$, then the corresponding values of $v$ can be recovered uniquely after cancelling $(1-x_0)u$ by $1-x_1$ on the left. This can be done uniquely since the group ring $R=K[F]$ has no zero divisors. So we are going to describe the values of $u$ in the equation. Obviously, they form a right $R$-module.

\begin{thm}
\label{descr}
Let $R=K[F]$ be a group ring of $F$ over a field $K$. Let 
$$
u_0=(1+x_0-x_1)(1-x_3),
$$
$$
u_1=(1-x_1)\phi(u_0)=(1-x_1)(1+x_1-x_2)(1-x_4),
$$
$$
\dots
$$
$$
u_k=(1-x_1)\phi(u_{k-1})=(1-x_1)\ldots(1-x_k)(1+x_k-x_{k+1})(1-x_{k+3})
$$
$$
\dots
$$

Then for any solution of the equation $(1-x_0)u=(1-x_1)v$ in $R$, the element $u$ belongs to the right $R$-module generated by $u_0$, $u_1$, $u_2$, \dots\,. Moreover, for any solution in the monoid ring $K[M]$, one can express it as $u=u_0r_0+u_1r_1+\cdots+u_kr_k$ for some $k$ and elements $r_i\in K[M_i]$ $(0\le i\le k)$.
\end{thm}

This immediately follows by induction on $\deg v$ from Lemma~\ref{summ} and formulas~(\ref{uuvv}). From this we have the description of solutions in the monoid ring $K[M]$. All solutions in $K[F]$ are obtained from them via a right multiplication by an element in $F$ according to Lemma~\ref{gig}.
\vspace{1ex}

Notice that to any relation between group generators $x_0$, $x_1$ we can assign a solution of the equation $(1-x_0)u=(1-x_1)v$ in the same way as we did it in the beginning of this Section. This is done uniquely up to right multiplication by an element of the group. So we can conclude that our description of solutions from Theorem~\ref{descr} gives some control over the set of group relations.
\vspace{1ex}

We finish this Section adding one more related result. Let $G$ be a group generated by a finite set $A=\{a_1,\ldots,a_m\}$. For any field $K$ we denote by $K[[G]]$ the space of infinite sums of the form
$$
\sum\limits_{g\in G}\alpha(g)\cdot g,
$$
where $\alpha(g)\in K$ are coefficients $(g\in G)$. Obviously, $G$ acts on the left and on the right on this space. These actions can be naturally extended to the group ring. This gives $K[[G]]$ the structure of a $K[G]$-bimodule. For the rest of this Section, we  assume that $K=\mathbb R$ will be the field of reals.
\vspace{1ex}

Let $\Gamma={\cal C}(G,A)$ be the left Cayley graph of $G$ with generating set $A$. Recall that a {\em flow} on a locally finite graph $\Gamma$ is a real-valued function $f\colon E\to\mathbb R$ such that $f(e^{-1})=-f(e)$. We say that $f(e)$ is the flow {\em through} the directed edge $e$. Given a vertex $v$, we define an {\em inflow} to it as a sum of flows through all edges with $v$ as a terminate vertex.

In~\cite{Gu22a} we formulated the following criterion of non-amenability in terms of flows. It can be deduced from other known criteria.

\begin{prop}
\label{critfl}
Let $G$ be a group with finite generating set $A$, and let $\Gamma={\cal C}(G;A)$ be its Cayley graph. The group $G$ is non-amenable if and only if there exist constants $C > 0$ and $\varepsilon > 0$, and a flow $f$ on $\Gamma$ with the following properties:
	
a$)$ The absolute value of the flow through each edge is bounded: $|f(e)|\le C$ for all $e\in E$;
	
b$)$ The inflow into each vertex is at least $\varepsilon$.	
\end{prop}

Slightly modifying the statement and the proof, we can assume that $\varepsilon=1$ and the inflow into each vertex equals exactly $1$. The key point here is to have uniformly bounded flows through all edges.

Having a flow on the Cayley graph with the above properties, we can form the sums $S_i=\sum\limits_{g\in G}k_i(g)g\in\mathbb R[[G]]$ for all $1\le i\le m$, where $k_i(g)$ denotes the flow through the directed edge labelled by $a_i$ from $a_ig$ to $g$ in the left Cayley graph $\Gamma$. It is easy to see that the sum $(1-a_1)S_1+\cdots+(1-a_m)S_m$ will be exactly the sum of all elements of the group: $\sum\limits_{g\in G}g$. This happens because all inflows into vertices are equal to $1$. The absolute values of the coefficients in the sums $S_i$ are uniformly bounded by a constant $C > 0$. Conversely, if $\sum\limits_{g\in G}g$ can be presented as a sum with bounded coefficients, then they bring a desired flow on the Cayley graph. So we have the following criterion in terms of group series.

\begin{thm}
\label{series}
Let $G$ be a group generated by $A=\{a_1,\ldots,a_m\}$. Then $G$ is non-amenable if and only if there exists an equality in the left $\mathbb R[G]$-module $\mathbb R[[G]]$ of the form
$$
\sum\limits_{g\in G}g=(1-a_1)S_1+\cdots+(1-a_m)S_m,
$$
where $S_1,\ldots,S_m$ are elements of $\mathbb R[[G]]$ with uniformly bounded coefficients.
\end{thm}

In partucular, this is applied to $F$, where we need an equality of the form $\sum\limits_{g\in F}g=(1-x_0)u+(1-x_1)v$, where $u$, $v$ are infinite group series from $\mathbb R[[F]]$ with uniformly bounded coefficients.

\end{document}